\documentclass[11pt]{amsart}
\usepackage{amsmath,amsthm,amssymb}

\usepackage{mathrsfs}
\let\mathcal\mathscr

\newcommand{\del}{\partial}

\newcommand{\cE}{{\mathcal E}}

\newcommand{\cO}{{\mathcal O}}

\newcommand{\cT}{{\mathcal T}}

\newcommand{\tA}{\widetilde{A}}

\newcommand{\tcT}{\widetilde{\mathcal T}}

\newcommand{\ra}{\rightarrow}
\newcommand{\lra}{\longrightarrow}

\newcommand{\G}{\Gamma}

\newcommand{\T}{\Theta}

\newcommand{\Sg}{\mathfrak{S}}

\def\Z{{\bf Z}}

\def\P{{\bf P}}

\def\C{{\bf C}}

\def\k{{\bf k}}

\def\div{\mathop{\rm div}\nolimits}
\def\Pic{\mathop{\rm Pic}\nolimits}

\def\Sym{\mathop{\rm Sym}\nolimits}
\def\Proj{\mathop{\rm Proj}\nolimits}

\newtheorem{theoremi}{Theorem}

\newtheorem{lemmai}[theoremi]{Lemma}
\newtheorem{propositioni}[theoremi]{Proposition}

\begin{document}

\title[Ampleness of intersections of translates of theta
divisors]{Ampleness of intersections of translates of theta divisors
in an abelian fourfold}

\author{O. Debarre}

\author{E. Izadi}

\address{Department of Mathematics\\
 Boyd
Graduate Studies Research Center\\
 University of Georgia\\
  Athens, GA
30602-7403\\ USA}

\email{izadi@math.uga.edu}

\address{IRMA -- Math\'ematique\\
 Universit\'e Louis Pasteur\\
 7 rue Ren\'e Descartes
 \\67084 Strasbourg Cedex\\ France}

\email{debarre@math.u-strasbg.fr}

\thanks{This work was done while O. Debarre was visiting the
Department of Mathematics of the University of Georgia, at the
invitation of E. Izadi. He is grateful to E. Izadi for her
hospitality. Both authors are grateful to the University of Georgia
for its support.}



\maketitle

\section*{Introduction}

Varieties  with ample cotangent bundle are interesting from many points of view. If $X$ is such a variety, defined over a field $\k$,
\begin{itemize}
\item (geometric) all   subvarieties of $X$ are of general type and there are only finitely many
rational maps from any fixed projective variety to $X$ 
(\cite{noguchisunada82});
\item (analytic) if $\k=\C$, any holomorphic map $\C\to X$ is constant (\cite{dem}, (3.1.));
\item (arithmetic) if $\k $ is a number field, the set of $\k$-rational points of $X$ is conjectured to be finite (see \cite{mo}; the analogous statement over function fields of curves is known to hold by \cite{noguchi82} or \cite{md}).
\end{itemize}
However, there
are relatively few concrete examples of such varieties. Bogomolov was the first to give a   general procedure to produce  such examples (his construction is explained in \cite{deamp}). In that article, more examples are constructed: it is shown that given 
 a principally polarized abelian variety $(A ,\T)$, an integer $n\ge\frac12\dim A   $, and  sufficiently ample (i.e.,
algebraically equivalent to sufficiently high multiples of $\T$)
general divisors $D_1, \dots, D_n$, the smooth
 variety $D_1\cap \cdots \cap D_n$   has ample cotangent bundle. In
this paper we prove an analogous result  for general abelian fourfolds. We work over an algebraically closed field $\k$.

\begin{theoremi}\label{thm1}
Let $(A ,\T)$ be a {\em general} principally polarized abelian
fourfold. For $a\in A$ general, the smooth surface $\T\cap\T_a$ has
ample cotangent bundle.
\end{theoremi}

Here $\T_a$ denotes the translate $\T+a$ of $\T$ by $a$. Our proof
shows that the conclusion of the theorem already holds on a general
Jacobian fourfold.

\section{The ampleness of $\T\cap\T_a$}

Let $(A ,\T)$ be a  principally polarized abelian fourfold. Assume
$\T\cap\T_a$ is a smooth surface. The cotangent bundle
$\Omega_{\T\cap\T_a }$ fits into the exact sequence of conormal and
cotangent bundles
\[
0\lra\cO_{\T\cap\T_a } (-\T )\oplus\cO_{\T\cap\T_a } (-\T_a
)\lra\Omega_A \vert_{\T\cap\T_a }\lra\Omega_{\T\cap\T_a }\lra 0
\]
Being a quotient of a trivial vector bundle, it is   generated by its
global sections, which are identified with $H^0 (\T\cap\T_a,\Omega_A
\vert_{\T\cap\T_a })\simeq H^0 (A,\Omega_A)$. To show that
$\Omega_{\T\cap\T_a }$ is ample, we must show that the associated
morphism
\begin{equation}\label{psia}
\psi_a:\P( \Omega_{\T\cap\T_a })\lra \P(H^0 (A,\Omega_A))
\end{equation}
is finite (we follow  Grothendieck's notation: given a coherent sheaf $\cE$ on a scheme, we set $\P(\cE)=\Proj(\Sym \cE)$). The set-theoretical inverse image of
a point of $\P(H^0 (A,\Omega_A))$, corresponding to a hyperplane in
$H^0 (A,\Omega_A)$, or to a line $\ell$ in $T_{A,0}$, is the set of
pairs $(\ell, x)$ where $x$ is in $\T\cap\T_a$ and $\ell$ is contained
in $T_{\T\cap\T_a ,x}$ (we identify $T_{A,0}$ and $T_{A,x}$   by
translation by $x$).

In other words, to prove that $\T\cap\T_a$ is ample, we must prove
that, for any nonzero $\del \in T_{A,0}$, the set of points
$x\in\T\cap\T_a$ such that $\del \in T_{\T\cap\T_a ,x}$ is finite.
We denote by $[\del]$ the point of $\P(H^0 (A,\Omega _A))$ determined by $\del$.

\section{The divisor $\T\cap\del\T$}

Let $(A ,\T)$ be a  principally polarized abelian variety and let
$\theta$ be a nonzero section of  $ \cO_A (\T ) $. We define, for any
$\del$ in  $T_{A,0}$, a section
$\del\theta$ of $\cO_{\T } (\T )$ by the requirement that for any
open set $U$ of $A$ and any trivialization
$\varphi:\cO_U\stackrel{{}_{{}_{\displaystyle\sim}}}{\lra}\cO_{\T }(\T
)\vert_U$, we have
$\partial \theta=\varphi(\partial(\varphi^{-1}(\theta)))\vert_\T$ in
$U\cap \T$. We
denote its zero locus by
$\T\cap\partial \T$. Set-theoretically, $\T\cap\del\T$ is the set of
points $x$ of $\T$ where the Zariski tangent space $T_{\T,x}$ contains
$\del$.

The differential of the isomorphism $A\to \Pic^0(A)$ induced by the
polarization $\T$ identifies $T_{A,0}$ with $T_{\Pic^0(A),0}\simeq
H^1(A,\cO_A)$.
The exact sequence
\[
0\lra\cO_A\lra\cO_A (\T )\lra\cO_{\T } (\T )\lra 0
\]
yields a composed isomorphism 
\begin{equation}\label{2}
   H^0(\T,\cO_{\T } (\T ))\stackrel{{}_{{}_{\displaystyle\sim}}}{\lra}
H^1(A,\cO_A)
\stackrel{{}_{{}_{\displaystyle\sim}}}{\lra}T_{A,0}
\end{equation}
whose inverse is given by 
$$\del  \longmapsto  \del\theta$$ 
 
 When $A$ has dimension $4$, the ampleness of the cotangent bundle of
$\T\cap\T_a$ is equivalent to  the following: {\em for {\em all}
nonzero $\del'\in T_{A,0}$, the scheme $\T\cap\del'\T\cap\T_a\cap\del'
\T_a$ is finite.}
 
 \medskip

For $\del\ne 0$, the scheme $\T\cap\del\T$ is a limit of intersections
of translates of $\T$ in the following sense. Let $m: A\times \T\ra A$
be the morphism $(x,y)\mapsto y-x$ and let $\cT=m^{-1}(\T)$. The
second projection $ A\times \T\to \T$ identifies the fiber at $a$ of
the first projection
$$\cT \lra
 A
$$
with $\T\cap\T_a$. If $\tA\to A$ is the blow-up of $0$, with
exceptional divisor $\P (\Omega_{A,0}) $, and $\tcT$ is the strict
transform of $\cT$ in $\tA\times\T\to A\times\T$, we obtain a  
family
$$ \tcT \lra \tA
$$
whose fiber at $[\del]\in\P (\Omega_{A,0})$ is isomorphic to
$\T\cap\del\T$. If $\T$ is irreducible, this is a flat family of
subvarieties of codimension $2$ of $A$.

We will study the ampleness of the cotangent bundle of $\T \cap\T_a$
by letting it specialize to $\T\cap\del\T$.

\section{The finiteness of $\T\cap\del\T\cap\del'\T\cap\del\del'\T$}

Let again $(A,\T)$ be a principally polarized abelian fourfold and let $\del$ be a nonzero element of $T_{A,0}$. If $\del'$ is in $T_{A,0}$,
we may define as above a section $\del\del'\theta$ of
$\cO_{\T\cap\del\T\cap\del'\T } (\T )$ whose zero locus we denote by
$\T\cap\del\T\cap\del'\T\cap\del\del'\T$.

The ampleness of the cotangent bundle of $\T\cap\del\T $ is equivalent
to the finiteness of the morphism
\[
\psi_\del:\P( \Omega_{\T\cap\del\T  })\lra \P(H^0 (A,\Omega
 _A))
\]
analogous to (\ref{psia}).  With our definitions, this means that {\em
for {\em all} nonzero $\del'\in T_{A,0}$, the scheme
$\T\cap\del\T\cap\del'\T\cap\del\del'\T$ is finite.} Unfortunately,
this never holds, because
this scheme has codimension at most $3$ for $\del'=\del$. We will
prove that for $A$ a general Jacobian and $\del$ general in $T_{A,0}$,
this is the only positive-dimensional fiber of $\psi_\del$.

Let $C$ be a smooth  curve of genus $4$, take $A = \Pic^3 C$,  and
let $\T\subset A$ be Riemann's theta divisor parametrizing effective
divisors of degree $3$ on $C$.

\begin{propositioni}\label{pp2}
For $C$ and $\del$ general, all fibers of the morphism
$\psi_\del$
are finite, except for that of $[\del]$.
\end{propositioni}

\begin{proof} Assume that $C$ is not hyperelliptic and identify
it with its canonical model in $\P^3=\P(H^0 (C,\omega_C )) =\P(H^0
(A,\Omega
_A))$, where it is the intersection of a quadric $Q$
 (which will be assumed to be smooth) 
and a cubic.  

Take $\del\in T_{A,0}$ nonzero such that $[\del]\notin C$. If, in addition, $C$ is not bielliptic, it can be
proved as in \cite{deb}, Th\'eor\`eme, that  
$S_\del=\T\cap\del\T$ is an integral, i.e., irreducible and reduced,
surface. Noting that the restriction $H^1(A,\cO_A)\to H^1(\T,\cO_\T)$
is bijective and using the isomorphism (\ref{2}), we obtain from the
exact sequence
\[
0\lra\cO_{\T }\stackrel{\del\theta }{\lra }\cO_{\T } (\T
)\lra\cO_{S_\del}(\T)\lra 0
\]
an exact sequence
$$\begin{matrix} 
0&\lra&\k&\stackrel{{}\cdot\del }{\lra } &T_{A,0}&\lra& H^0 (S_\del ,\cO_{S_\del} (\T ))\\
&&&&\del'&\longmapsto&\del'\theta
\end{matrix}
$$
{\em Assume $\del$ and $\del'$ are linearly independent.} The section
$\del'\theta$ is then  nonzero on $S_\del$, and its zero locus
$\Gamma =\T\cap\del\T\cap\del'\T$ is a curve.  Similarly, since
$H^1(A,\cO_A)\to H^1({S_\del},\cO_{S_\del})$ is bijective, we obtain
from the exact sequence
\[
0 \lra \cO_{S_\del}\stackrel{\del'\theta }{\lra }\cO_{S_\del} (\T
)\lra\cO_{\G}(\T)\lra 0
\]
  a coboundary map
$ 
H^0(\G,\cO_{\G}(\T))\to T_{A,0}$ that sends
$\del\del' \theta$ to $\del$. This section is in particular nonzero. This shows that, {\em if $\G$ is integral,  
$\T\cap\del\T\cap\del'\T\cap\del\del'\T$ is finite,} which is what we need to prove.

The projectivization of the tangent space to $\T$ at a point
corresponding to a divisor $D$ of degree $3$ such that $h^0 (C,D) = 1$
is the plane spanned in $\P^3$ by the points of $D$.
The underlying
reduced curve $\Gamma_{\rm red}$ therefore parametrizes effective divisors  
of degree $3$ on $C$ that lie in a plane that contains the line
$\ell_{\del,\del'} $ spanned by $[\del]$
and $[\del']$. We will distinguish several cases depending on the
relative positions of $[\del]$, $[\del']$, and $C$ in $\P^3$.

We first introduce some notation, following \cite{I13}:
given a pencil $g^1_e$ on $C$ with reduced base locus, we define, for
any $d\in\{1,\dots,e\}$, a reduced curve in the $d$-th symmetric
power $C^{(d)}$ by setting
$$X_d(g^1_e)=\{p_1+ \cdots+ p_d\in C^{(d)}\mid\exists D\in C^{ (e-d)}\ D+ p_1 +\cdots + p_d\in g^1_e\}.
$$

\subsection{Case 1: $\ell_{\del,\del'}\cap C =\varnothing$}

The planes containing $\ell_{\del,\del'}$ cut on $C$ the
divisors of a base-point-free $g^1_6$  contained in $|\omega_C | $, and the curve
$\Gamma_{\rm red}$ is the image in $\T$ of the curve
$X_3(g^1_6)\subset C^{(3)}$. It follows from  \cite{acgh}, Lemma
VIII.(3.2)  that the cohomology class of $\Gamma_{\rm
red}$ is   $[\T]^3$, so $\Gamma$ is reduced.

The associated map $\phi : C\to (g^1_6)^*=\P^1$ coincides with the
projection of $C\subset\P^3$ from the line $\ell_{\del,\del'}$. The lemma below
shows that the monodromy group $G$ of $\phi$ is the full symmetric
group $\Sg_6$. It implies that $\Gamma$ is integral, and we are done in this case.

\begin{lemmai}
 For $C$ general
  and $[\del]\notin Q$, the group $G$   is $2$-transitive and contains
a simple transposition.
\end{lemmai}

\begin{proof}
The $2$-transitivity of $G$ is equivalent to the irreducibility of the curve
$
X_2 (g^1_6 ) $ in $ C^{ (2) }
$.

Let $\pi:C^2\to C^{ (2) }$ be the quotient map. For any divisor
(resp. divisor class) $D$ on $C$, let $C_D$ be the unique divisor
(resp. divisor class) on $C^{ (2) }$ such that $\pi^*C_D=p_1^*D+
p_2^*D$. Let $\delta$ be the unique divisor class on $C^{ (2) }$ such
that $\pi^*\delta$ is linearly equivalent to the diagonal of
$C^2$. We have the linear equivalence $X_2 (g^1_6)\equiv
C_{ g^1_6} -\delta$ (\cite{I13}, Lemma 2.1). Moreover, $\delta^3=-3$ and $\delta\cdot C_D=\deg (D)$.

Assume $C$ is sufficiently general so  that the map
$$\begin{matrix} 
\Pic(C)\oplus\Z&\lra&\Pic(C^{ (2) })\\
(D,b) &\longmapsto&C_D-b\delta
\end{matrix}
$$
is bijective.
If $X_2 (g^1_6 )$ is reducible,  write the divisor class of a nontrivial   union of  components, say $Y$, as $C_D-b\delta$, so that the class of $X_2 (g^1_6 ) -Y$ is $C_{g^1_6 -D}- (1-b )\delta$. Replacing $Y$ with $X_2 (g^1_6 )-Y$ if necessary, we may assume $b\ge 0$.

We now use  \cite{I13},  Appendix 6.1: for any divisor $E$ on $C$, we have
$$ 
H^0(C^{ (2) }, C_E )\simeq \Sym^2H^0(C,E)
\quad{\rm and}\quad H^0(C^{ (2) },
C_E-\delta )\simeq \bigwedge^2H^0(C,E)$$
It follows that if $E$ is effective and $h^0(C,E)= 1$, the linear system $|C_E- \delta|$  is empty, and  $|C_E|=\{C_E\}$. Since our $g^1_6$ has no base points, $X_2 (g^1_6 )$ contains no such curve. It follows that $D$ moves in a pencil, hence $\deg(D)\ge 3$ since $C$ is not hyperelliptic. Since the diagonal is not a component of $X_2 (g^1_6 ) $, we must have $(C_{g^1_6 -D}- (1-b )\delta)\cdot \delta\ge0$, hence $3b\le 9-\deg (D)$.

If $\deg(D)\ge 4$, we get $b\le 1$ but this is absurd since $|C_{g^1_6 -D}- (1-b )\delta|$ is then empty.  Hence $D$ is a $g^1_3$ and $b\le 2$. By \cite{I13},  Appendix 6.3, the vector subspace $ H^0(C^{ (2) }, C_{g^1_3}-2\delta )\subset H^0(C^{ (2) }, C_{g^1_3} )$ is isomorphic to the space of quadratic forms that vanish on the image of $C\to({g^1_3})^*$, hence vanishes. We get   $b\in \{0,1\}$ and, replacing $Y$ with $X_2 (g^1_6 )-Y$ if necessary, 
 $ Y \equiv C_{g^1_3}-\delta$. More precisely, $Y=X_2 (g^1_3)$.
  The $ g^1_3$ is given by one of the rulings of the quadric $Q$, hence $ X_2 (g^1_3)$ may be contained in  $ X_2 (g^1_6)$ only if the line $\ell_{\del,\del'}$ meets all lines of this ruling. Since  $[\del]$ is not in $Q$, this cannot happen and $X_2 (g^1_6 )$ is irreducible.

To prove that $G$ contains a simple transposition, we check that
for $C$ general, there is a  point $p\in C$ such that $\phi : C\to
\P^1$ ramifies simply at $p$ and $p$ is the only ramification point of
$\phi$ in its fiber.

 The degree of the ramification locus is $6 + 6\cdot 2 = 18$. If all
the ramification points are either nonsimple or their fiber contains
other ramification points, the support of the branch locus of $\phi$
in $\P^1  $ contains at most $9$ points. Such $6$-fold covers of
$\P^1$ depend on at most $9-3 = 6$ parameters. Therefore, for a
sufficiently general choice of $C$, the map $\phi$ will have at least
$3$ ramification points with the desired property.
\end{proof}

  We assume from now on that $[\del]$  lies on no trisecants
($[\del]\notin Q$) or tangents to $C$.

\subsection{Case 2: $\ell_{\del,\del'}\cap C =\{ p\}$}

Here we mean that the line $\ell_{\del,\del'}$ is {\em not} tangent to
$C$. In this case we have
\[
\Gamma_{\rm red } = X_3 (g^1_5)\cup \bigl( X_2 (g^1_5 )+p\bigr)\subset C^{ (3) }
\]
where $g^1_5\subset |\omega_C |$ is the base-point-free pencil cut
on $C$ by planes through $\ell_{\del,\del'}$. As before, we see that
$\Gamma$ is reduced. The involution
\[
\tau: x+y+z\longmapsto K_C-x-y-z
\]
exchanges $X_3 (g^1_5 )$ and $X_2 (g^1_5 )+p$. A similar (simpler)
calculation as before shows that $X_2 (g^1_5)$, hence also $X_3 (g^1_5
)$, is irreducible. As the scheme $\Gamma\cap\del\del'\T$ is invariant
under $\tau$, we see that if it contains one component of $\Gamma$, it
also contains the other. This is therefore not possible, hence this scheme is finite.

\subsection{Case 3: $\ell_{\del,\del'}\cap C =\{ p,q\}$}

Here we mean that the line $\ell_{\del,\del'}$ intersects $C$ in
exactly two distinct points $p$ and $q$.

Let $\del_p $ and $\del_q$ be elements of $T_{A,0}$ mapping to $p$ and
$ q $ respectively.
Let $W_p$ be the image in $\T$ of $p+C^{(2)}\subset C^{ (3) }$. We have
\[
\T\cap\del_p\T = W_p\cup (K_C-W_p) = W_p\cup\tau (W_p)
\]
Since $[\del]\notin Q$, the linear system $|K_C-p-q|$ is a
base-point-free $g^1_4$ and the curve $X_2 (K_C-p-q)$ is irreducible
as before. We have
$\Gamma =\T\cap\del\T\cap\del'\T=\T\cap\del_p\T\cap\del_q\T $ and we check that this curve is
reduced and has four irreducible components:
\[
\begin{array}{rclcrcl}
\G_1&=&p+q+C&&\G_2&=&p+ X_2 (K_C-p-q)\\
\tau(\G_1)&=&X_3 (K_C-p-q)&&\tau(\G_2)&=&q+ X_2 (K_C-p-q)
\end{array}
\]
If $\T\cap\del\T\cap\del'\T\cap\del\del'\T$ contains a component of
$\Gamma$, it also contains its image by $\tau$.
 It will therefore be enough for our purpose to show that  the section
$ \del\del'\theta$ of $\cO_\G(\T)$  vanishes identically neither on
$\G_1$, nor on $\G_2$ (both contained in $W_p $).

Let $\iota_{p+q}$ be the embedding $x\mapsto p+q+x$ of $C$ into
$A$, with image $\G_1$. We have $\iota_{p+q}^*\Theta\equiv  K_C-p-q$. Let $p+p_1+p_2+p_3$ and
 $q+q_1+q_2+q_3$ be the divisors of $|K-p-q |$ containing $p$ and
$q$. For a sufficiently general choice of $\del$, these two divisors
will be reduced and disjoint.

 \begin{lemmai}
 The section $ \iota_{p+q}^*\del_p\del_q\theta$   vanishes identically   and 
\begin{eqnarray*}
 \div(\iota_{p+q}^*\del_p^2\theta)&=&p+p_1+p_2+p_3\\
\div(\iota_{p+q}^* \del_q^2\theta)&=&q+q_1+q_2+q_3
\end{eqnarray*}
\end{lemmai}

\begin{proof}Let $\lambda\in\k$ and set $\del_\lambda=\lambda\del_p
+\del_q $. The support of
$$\div
(\del_p\del_\lambda\theta)=\T\cap\del_p\T\cap\del_q\T\cap\del_p\del_\lambda\T
=\T\cap\del_p\T\cap\del_\lambda\T\cap\del_p\del_\lambda\T
$$ is the set of points of $\T\cap\del_p\T = W_p\cup (K_C-W_p)$ whose
tangent space contains $\del_\lambda$.

 It contains $p+q+x$ if the line $\langle q,x\rangle$ contains
 $[\del_\lambda]$. In particular, $ \del_p\del_q\theta$ vanishes
 identically on $p+q+C$ and $ \del_p \del_\lambda\theta(2p+q)=0$ for
 all $\lambda$. This implies $ \del_p^2\theta(2p+q) =0$. Moreover, $
 \del_p
\del_\lambda\theta(p+2q)\ne 0$ if $\lambda\ne0$. In particular,
$\iota_{p+q}^*\del_p^2\theta$ is a nonzero section of $\cO_C(K_C-p-q)$
that vanishes at $p$, hence  the lemma.
  \end{proof}

Write $\del =\lambda\del_p +\mu\del_q$ and $\del' =\lambda'\del_p
+\mu'\del_q$, so that
\[
\del\del'\theta =  \lambda\lambda'\del_p^2\theta + (\lambda\mu'
+\lambda'\mu )\del_p\del_q\theta +\mu\mu'\del_q^2\theta.
\]
Since $[\del]$ is not on $C$, both $\lambda$ and $\mu$ are not zero,
hence   $\del\del'\theta$ {\em does not vanish identically on $\G_1$.}
We have
\begin{eqnarray*}
\G_1\cap\G_2&=&\{p+q+q_1,p+q+q_2,p+q+q_3\}\\
\G_1\cap\tau(\G_2)&=&\{p+q+p_1,p+q+p_2,p+q+p_3\}\\
\tau(\G_1)\cap\G_2&=&\{\tau(p+q+p_1),\tau(p+q+p_2),\tau(p+q+p_3
)\}\end{eqnarray*}
The section $\del_p\del_q\theta $ does not vanish identically on
$\Gamma$, hence does not vanish identically on $\Gamma_2$.
At $p+q+q_1$, both $\del_p\del_q\theta $ and $\del_q^2\theta $ vanish,
but $\del_p^2\theta $ does not. At $\tau(p+q+p_1)$, both
$\del_p\del_q\theta $ and $\del_p^2\theta $ vanish, but
$\del_q^2\theta $ does not. It follows that the sections
$\del_p^2\theta\vert_{\G_2}$, $\del_p\del_q\theta\vert_{\G_2}$, and
$\del_q^2\theta\vert_{\G_2}$ are linearly independent, hence
$\del\del'\theta$ does not vanish identically on $ \G_2$.

We have proved that in all cases, the zero set of $\del\del'\theta$ on $\G$ is
finite. This completes the proof of Proposition \ref{pp2}.
\end{proof}

\section{The scheme $\T\cap\del\T\cap\del^2\T$}

The fiber of $\psi_\del^{-1}([\del])$  is one-dimensional, equal to
$
\T\cap\del\T\cap\del^2\T
$. We now study this curve. Let $p$ be a general point of $C$. As
above, we see that  $\T\cap\del_p\T\cap\del_p^2\T $
has three irreducible components whose reduced underlying curves are
$$
\begin{array}{rcl}
\G_1=2p+C\hskip-5mm&&\hskip-5mm\tau(\G_1)=X_3 (K_C-2p)\\
&\hskip-6mm\G_2=\tau(\G_2)=p+ X_2 (K_C-2p)\hskip-15mm
\end{array}
$$
and $\G_1$ and $\tau(\G_1)$ have multiplicity $1$, whereas $\G_2$ has
multiplicity $2$.

\begin{lemmai}\label{l4}
For $\del$  general, $\T\cap\del\T\cap\del^2\T$ contains no translates of $C$.
\end{lemmai}

\begin{proof}
A translate of $C$ is contained in $\Theta$ if and only if it is of the type $x+y+C$, with $x$, $y\in C$. It is contained in $\T\cap\del\T$ if and only if for
every $t\in C$, the plane $\langle x,y,t\rangle$ contains
$[\del]$. This is only possible if $[\del]$ is on the line $\langle
x,y\rangle$. For $\del$ general, there are exactly six distinct
secants to $C$ that contain $[\del]$, none of which is trisecant or
tangent. So there are exactly six distinct translates, say $x_i + y_i
+C$, for $i\in\{1,\dots,6\}$, contained in $\T\cap\del\T$. Since the
set of secants to $C$ is irreducible, if one of these translates is
contained in $\T\cap\del\T\cap\del^2\T$ for $\del$ general, they all are. This implies
\[
\T\cap\del\T\cap\del^2\T =\bigcup_{ i=1 }^6 (x_i + y_i +C)
\]
which is not possible since a general $\T\cap\del\T\cap\del^2\T$ has
at most four irreducible components by the description of
$\T\cap\del_p\T\cap\del_p^2\T$ above.
\end{proof}

Since $4p$ is not contained in a plane in $\P^3$, the curves
$\G_1$ and $\tau(\G_1)$ defined above are disjoint. Therefore, it follows from Lemma
\ref{l4} that if a general $\T\cap\del\T\cap\del^2\T$ is nonintegral,
it is of the form $\Gamma_0\cup\tau (\Gamma_0 )$, where $\Gamma_0$ is
integral, with cohomology class $\frac12 [\T ]^3$, and distinct from
$\tau(\G_0)$.

\section{Proof of Theorem \ref{thm1}}

We keep the same assumptions and notation as before. Let $a$ be
general in $A=\Pic^3(C)$. If for some nonzero $\del'$, the scheme
$\T\cap\T_a\cap\del'\T\cap\del'\T_a$ has dimension $1$, it contains a
curve $\G_a$ that is stable by the involution $x\mapsto a-x$. When $a$
specializes to a general $[\del]$, this involution specializes to
$\tau$ and $\G_a$ must specialize as a set to $\Gamma_0\cup\tau
(\Gamma_0 )$. Since this curve has the same cohomology class as the
curve $\T\cap\T_a\cap\del'\T$, this means that the section
$\del'\theta_a$ vanishes identically on the curve
$\T\cap\T_a\cap\del'\T$ and this is absurd.

It follows that $\psi_a$ is finite, hence the cotangent bundle of
$\T\cap\T_a$ is ample.



\end{document}